\documentclass[12pt]{article}
\usepackage{setspace}
\usepackage{graphicx,color}
\usepackage{amsmath}
\usepackage{amssymb}
\usepackage{amsthm}
\usepackage{epsfig}
\usepackage{latexsym}      
\newtheorem{theorem}{Theorem}
\newtheorem{lemma}[theorem]{Lemma}

\theoremstyle{definition}
\newtheorem{definition}{Definition}
                          
\newtheorem{conjecture}{Conjecture}
\newtheorem{question}{Question}

\newcommand{\R}[1]{\ensuremath{\mathbb{R}^{#1}}}
\newcommand{\C}[1]{\ensuremath{\mathbb{C}^{#1}}}
\newcommand{\D}[1]{\ensuremath{\mathbf{D}_{#1}}}
\newcommand{\K}[1]{\ensuremath{\mathbf{k}_{#1}}}
\newcommand{\LL}{\ensuremath{\widetilde{\K{p/q}}}}
\newcommand{\LM}{\ensuremath{\widetilde{\K{m}}}}
\newcommand{\N}{\ensuremath{\mathbf{n}_{m}}}
\newcommand{\RGC}{\ensuremath{\lbrace (z,w) : z \in \mathbb{R}, w \in \mathbb{R}, z^2 + w^2 =1 \rbrace}}

\date{}
\title{Great circle links and virtually fibered knots}
\author{Genevieve S. Walsh \thanks{University of Texas at Austin, Department of Mathematics. email: gwalsh@math.utexas.edu}}
\begin{document}  
\maketitle
 
\abstract{We show that all two-bridge knot and link complements are virtually fibered. We also show that spherical Montesinos knot and link complements are virtually fibered.  This is accomplished by showing that such manifolds are finitely covered by  great circle link complements. }

\section{Introduction}

A three-manifold $M$ is  {\it fibered} if it admits a submersion $p: M \rightarrow S^1$.  Equivalently, $$M = \frac{F \times I}{(x,0) \sim (\psi(x),1)}$$ 
 where $F$ is a surface and $\psi$ is an orientation preserving isomorphism from $ S$ to itself.  Fibered  three-manifolds are very well understood, since they are completely determined by the surface $F$ and the monodromy $\psi$.  For example, if the monodromy of the fibration is pseudo-Anosov, then the fiber bundle admits a hyperbolic structure, \cite{ThurstonHypII}.  
Although fibered manifolds are particularly nice, there are of course many manifolds that are not fibered.  A manifold $M$ is {\it virtually fibered} if there is a finite cover of $M$ that is fibered.  This includes the case that $M$ is itself fibered.   It is conjectured that many manifolds that are not fibered are virtually fibered, \cite{ThurstonBulletin}.  The "virtually fibered" conjecture of W. P. Thurston is as follows:

\begin{conjecture} \label{virtfiberconj} \it{ Let $M$ be a compact orientable irreducible three-manifold whose fundamental group is infinite and contains no non-peripheral $Z \times Z$ subgroup. Suppose that any boundary components of $M$ are tori.  Then $M$ is virtually fibered. }

\end{conjecture}

For example, the conjecture would imply that any hyperbolic knot or link complement is either fibered or virtually fibered. 
In contrast,  Seifert fibered spaces are virtually fibered exactly when either the  Euler number of the fibration is zero or the orbifold Euler characteristic of the base orbifold is zero. See \cite{Gabaifibered} for a proof.  Since a finite cover of a fibered manifold is also fibered, the property of being virtually fibered is an invariant of the commensurability class.    Although there is not a lot of evidence for Conjecture \ref{virtfiberconj}, there are examples of hyperbolic three-manifolds that are known to be virtually fibered but not fibered.  The first non-trivial example was a link complement given by  Gabai in \cite{Gabaifibered}. Examples of  non-Haken virtually fibered manifolds are given in \cite{Reidfibered}.  The first examples of non-fibered virtually fibered knot complements were given by  Leininger, in  \cite{Chrisfibered}.  Here we exhibit a large class of virtually fibered knot and link complements, namely two-bridge knots and links and spherical Montesinos knots and links.   This provides  evidence for Conjecture \ref{virtfiberconj}.  Our main theorem is as follows. 
\vskip .1 in
{\bf Theorem \ref{twobridgevirtual}} {\it Every two-bridge knot complement and non-trivial two-bridge link complement is virtually fibered.}

\vskip .1 in
This is proven in section \ref{two-bridge}, by showing that every such knot or link complement is finitely covered by a great circle link complement.  Theorem \ref{twobridgevirtual} implies that the knot complements in \cite{Chrisfibered} are virtually fibered, but the approach is completely different.   We have a similar theorem for spherical Montesinos links. A {\it spherical Montesinos link} is a Montesinos link $k$ where the double branched cover of $S^3$ branched along $k$ is a spherical Seifert fibered space.  Note that a two-fold cover of $S^3$ branched along a knot or link as in Theorem \ref{twobridgevirtual} is a lens space, which is also spherical. 

\vskip .1 in 
{\bf Theorem \ref{montesinosfiber}} {\it The complement of any spherical Montesinos knot or link is virtually fibered.}

\vskip .1 in

This is proven in section \ref{montesinossection}.  By using, for example, the Alexander polynomial and computations of Hatcher and Oertel in \cite{HO}, we can determine which knots through nine crossings are fibered and which knots are spherical Montesinos or two-bridge. For example, the knot $9_{25}$ is hyperbolic, spherical Montesinos,  and not fibered. 
We conclude that, with five possible exceptions, all knots though nine crossings  are either fibered or virtually fibered.  The exceptional cases not covered by the results in this paper are $\lbrace 9_{38}, \, 9_{39}, \, 9_{41}, \, 9_{46}, \, 9_{49} \rbrace$.  There are many 10 crossing knots not covered by the results of this paper.  

 In section \ref{questions} we discuss some related questions, in particular, if these methods apply to other knot complements. 

\vskip .2 in

\section{Great circle links}\label{GCsection}

We begin with a definition.  

\begin{definition} A link in $S^3$ is a {\it great circle link} if all of the components are geodesics in $S^3$.
\end{definition} 
 We will consider $S^3$ as the the unit three-sphere in $\R{4} = \C{2}$.  Then every great circle is the intersection of a 2-plane through the origin in $\R{4}$  and the three-sphere.  Note that great circles are invariant under the transformation $x \mapsto -x$ from $S^3$ to itself.  Therefore, configurations of great circles are the same as configurations of lines in $\mathbb{RP}^3$. These configurations have been studied extensively, in particular by the Russian school, and are completely classified up to seven components, using combinatorial methods in \cite{RP3} \cite{sevenlines}.  Great circle links were later directly classified up to five components, using the geometric structures on their complements in \cite{Gthesis}.   The simplest examples of great circle links are Hopf links, links where every component is a fiber of the same Hopf fibration. However, there are many examples of great circle links whose complements admit hyperbolic structures, in particular many of the examples given here.   For the purposes of this paper, the most relevant feature of great circle links is the following: 

\begin{theorem} 
\label{fibered} All great circle link complements are fibered. 
\end{theorem}

\begin{proof} Let $c$ be any component of a great circle link $L$.  
Then $c$ is the intersection of $S^3$ and a two-plane $P_c$ in $\R{4}$ through the origin.   Then there is an  isometry of $\R{4}$ that takes $P_c$ to the plane defined by the vectors $(
1, \,0, \,0, \,0 )$ and $( 0,\, 1,\,0,\,0) $ in $\R{4}$.  This isometry takes $L$ to an equivalent great circle link. 
 Therefore, we may assume that $c$ is the great circle $\lbrace (x,\, y, \, 0, \, 0): x^2 + y^2 =1 \rbrace$ in $S^3 \subset \R{4}$.    The complement of $c$ in $S^3$ is an open solid torus that is fibered by the half-planes  $H_\theta = \lbrace (x,y, r \cos \theta, r \sin \theta ) \in S^3 , r >0 \rbrace $.   These are hemispheres of great spheres and are totally geodesic. Therefore if a geodesic arc intersects $H_\theta$ in a point it must do so transversely.

We claim that any component of $L -c $ intersects each $H_\theta$ exactly once.  Assume that some component does not intersect one of the hemispheres.  Then it does not link $c$.  But every two great circles in $S^3$ link with linking number $\pm 1$.  To see this, note that the union of the bases for the two corresponding two-planes through the origin forms a basis for $\R{4}$.  Then there is a linear transformation of $\R{4}$ that takes this basis to the standard one, and takes the two great circles to $\lbrace (a, \,b, \,0, \,0): a^2+b^2 =1 \rbrace$ and  $\lbrace (0, \,0, \,d, \,e): d^2+e^2 =1 \rbrace$.  These two circles link with linking number $ \pm 1$, depending on their orientations.  Therefore, every component of $L-c$ must intersect each $H_\theta$ at least once.

If two points $\bar x_1 = (x_1,y_1,r_1 \cos \theta,r_1 \sin \theta )$ and $ \bar x_2 = (x_2,y_2,r_2 \cos \theta,r_2 \sin \theta )$ in $S^3$ are in one of these hemispheres, then they are not antipodal. Since they are linearly independent points in $\R{4}$, the shortest geodesic between the two points is the intersection of $S^3$ and the plane in $\R{4}$ spanned by $\bar x_1$ and $\bar x_2$.  This geodesic is contained in the sphere $H_\theta \cup H_{-\theta} \cup c$ and intersects  $c$ in two antipodal points.   Therefore a geodesic in the complement of $c$ must intersect any $H_\theta$ exactly once.

Thus the map $(S^3 -L) \rightarrow S^1$ which sends a point in $H_{\theta}$ to $\theta$ is a fibration. The fibers are $(p-1)$-punctured disks, where $p$ is the number of components in $L$.  There are $p$ different such fibrations for a given link, one for each component. 
 
  \end{proof}

\section{Two-bridge knot and link complements \label{two-bridge}} 

\begin{definition}

Let $B$ and $B'$ be two balls, each with two standard interior vertical arcs marked. Let $h_{p/q}$ be the homeomorphism of their boundaries that takes a vertical arc (slope $\infty$) on $\partial B$ to an arc of slope $p/q$  on $\partial B'$.
The  rational knot $\K{p/q}$ is the union of the marked arcs in $S^3$ obtained by gluing together  $B$ and $B'$ by $h_{p/q}: \partial B \rightarrow \partial B'$. 
\end{definition}

\begin{figure}[htb]
\centerline{\epsfig{file=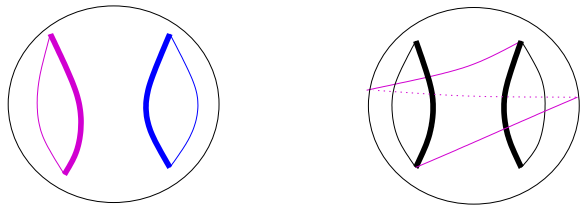, width=3 in}}
\caption{The homeomorphism $h_{1/3}$ takes arc of slope $\infty$  to an arc of slope 1/3. \label{fig:knota}}
\end{figure}

\begin{figure}[h]
\centerline{\epsfig{file=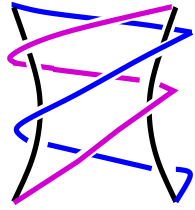, width=1 in}}
\caption{The knot $\K{1/3}$, the trefoil.\label{knotb}}
\end{figure}

The fact that every two-bridge knot can be put in this form was proven by Schubert in \cite{schubertknot}, which we record as 
\begin{theorem} 
Every two-bridge knot or link can be written as \K{p/q} for some $p/q$ with $ (p,q)=1$. 

\label{rationalform}
\end{theorem}

Note that the rational link with $p/q =1/0= \infty$ is the trivial two component link .  We want to exclude this case.

\begin{theorem} \label{twobridgevirtual} Every two-bridge knot complement and non-trivial two-bridge link complement is virtually fibered.
\end{theorem}

\noindent
Theorem \ref{twobridgevirtual} follows immediately from Theorem \ref{fibered} and the following. 

\begin{lemma} \label{GCLcover} Every two-bridge knot complement and non-trivial two-bridge link complement has a finite cover which is a great circle link complement. 
\end{lemma} 
\begin{proof}

Note that if $q$ is even, $\K{p/q}$ is a link, and if $q$ is odd, $\K{p/q}$ is a knot.  We are assuming that $q \neq 0$.  To prove Lemma \ref{GCLcover}, we make use of an alternative picture of great circle links in which we consider $S^3$ as the unit three-sphere in \C{2}.  Then we will refer to the great circles $\lbrace  (z, \,0), z \in \C{} \rbrace \cap S^3$ and $\lbrace (0, \, w), w \in \C{} \rbrace \cap S^3$ as the $z$-axis and the $w$-axis, respectively,  of $S^3$.  Consider the real great circle $g =$ \RGC. We can move $g$ around by the isometry   $$\phi_{p/q} : (z,w) \rightarrow (e^{\frac{2\pi i}{q}}z,  e^{\frac{2 \pi i p}{q}} w).$$ Since $\phi_{p/q}$ is an isometry of $S^3$, the image of the real great circle $g$ is another geodesic.  Let $p$ and $q$ be relatively prime.  Then the components of the orbit of $g$ under the action of $\phi_{p/q}$ either do not intersect or are identical.

When $q$ is odd, we let $\D{p/q}$ be the $q$-component great circle link that is the orbit of  $g$ under the isometry $\phi_{p/q}$.  Note that this  action on $S^3$ yields the lens space $L(q,p)$ as quotient.  We will denote this lens space as $L_{p/q}$.   This action leaves the $z$ and $w$ axes of $S^3$ invariant.  Since the real great circle intersects both of these axes in two points, and every component is a geodesic, the link is determined by the images of these intersection points.  Therefore, we can easily draw the link for small values of $p$ and $q$.  Figure \ref{D25picture} shows the link $\D{2/5}. $ 

\begin{figure}[htb]
\begin{center}
\input{5-2link.pstex_t}
\end{center} 
\caption{The link $\D{2/5}$
\label{D25picture}}
\end{figure}

Observe that  if $q$ is odd, the link  $\D{p/q}$ will intersect the $z$-axis in the points $(e^{\frac{k\pi i}{q}},0$ and the $w$-axis in the points $(0, e^{\frac{l \pi i}{q}})$, for $k,l \in \lbrace 1,2, ... 2q \rbrace$. We claim that the $\frac{q+1}{2}$th  component in the orbit of $\phi_{p/q}$ intersects the $z$-axis in $(e^{\frac{\pi i}{q}},0)$ and its antipodal point.  This is just because  $\frac{q+1}{2} (\frac{2 \pi}{q}) = \frac{\pi}{q} \mod \pi$.  This link component will intersect the $w$-axis in $(0, e^{\frac{p \pi i}{q}})$ and its antipodal point. Therefore, the $\frac{q+1}{2}$-th power of the isometry  $\phi_{p/q}$ has the effect of rotating the $z$-axis by $\pi/q$ and rotating the $w$-axis by $p \pi/q$.  We can draw the standard projection (where the $z$-axis is coming out of the page) as in figure \ref{D25standard}.  Start by drawing the $w$-axis (dotted) and the real great circle $g$.  Then the next component along the $z$-axis will be rotated $2 \pi / 5$ from $g$ along the $w$-axis.  We continue in this way to obtain figure \ref{D25standard}.

\begin{figure}[htb]
\begin{center}
\epsfig{file=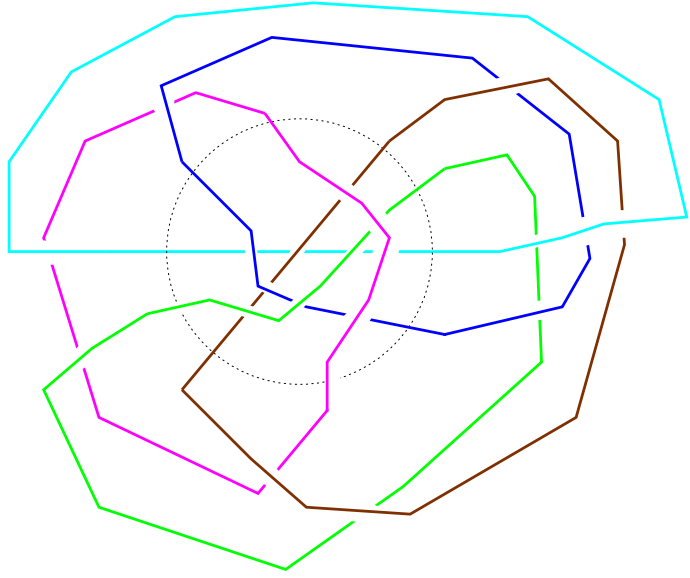,height=3truein}

\end{center} 
\caption{The standard projection of the link $\D{2/5}$
\label{D25standard}}
\end{figure}

Let $g_{a,b}$ be the great circle that intersects the the $z$-axis in the points $(\pm e^{ a i},0)$  and the $w$-axis in the points $(0, \pm e^{b i})$.  Then the real great circle $g$ is $ g_{0,0}$.
 In the case that $q$ is even, we define the link $\D{p/q}$ to be the orbit of the real great circle under $\phi_{p/q}$ union the orbit of $g_{\frac{\pi}{q},\frac{p \pi}{q}}$ under the action of $\phi_{p/q}$.  These two orbits will never intersect.  Like the case when $q$ is odd the resulting great circle link will intersect the $z$-axis in the points $(e^{\frac{k\pi i}{q}},0)$ and the $w$-axis in the points $(0, e^{\frac{ l \pi i}{q}})$ for $k,l \in \lbrace 1,2, ... 2q \rbrace$.  However, in this case, there will be two orbits under the action of $\phi_{p/q}$, one associated to the real great circle, where $k$ and $l$ are  always even, and one associated to the great circle $g_{\frac{\pi}{q},\frac{p \pi}{q}}$, where $k$ and $l$ are always odd. Note that in this case, component $x$ and component $x+ q/2$ under the orbit of $\phi_{p/q}$ will be the same component of $\D{p/q}$.  $\D{p/q}$ is again a $q$-component link.  
We will show that the link complement $S^3 - \D{p/q}$ covers a two-bridge knot or link complement.

First note that the two-fold branched cover of $S^3$ branched along a two-bridge knot or link $\K{p/q}$ is a lens space.  (see figure \ref{knotlens})  The two solid tori that are glued together so that a meridian curve is glued to a curve of slope $p/q$.  We will call this $L_{p/q}$. Schubert proved that  \K{p/q} and \K{p'/q'} are equivalent (as unoriented links) if and only if $q = q'$ and $p^{\pm1} = \pm p'(\mod q) $.  This is exactly when their associated lens spaces are homeomorphic. 
 
\begin{figure}[h]
\centerline{\epsfig{file=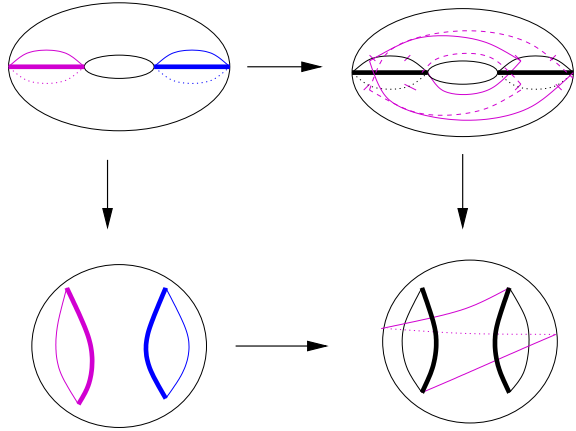, width=3 in}}
\caption{The two-fold cover of $S^3$ branched along the knot or link $K_{p/q}$ is $L_{p/q}$.\label{knotlens}
}
\end{figure}

Now let  $\LL$  be the pre-image of the branching locus in $L_{p/q}$.  Then the two-bridge link complement $S^3 -\K{p/q}$ is covered by $L_{p/q} - \LL$.   

We claim that the great circle link complement $S^3-\D{p/q}$ covers $L_{p/q} - \LL$. To prove this we will show $(S^3, \D{p/q})$ covers $(L_{p/q}, \LL)$ as a map of pairs.  Consider the action of $\phi_{p/q}:S^3 \rightarrow S^3$ where $\phi_{p/q}(z,w)= (e^{\frac{2\pi i}{q}} z ,e^{\frac{2\pi i p}{q}} w)$. A fundamental domain for this action is the union of the two wedges $\lbrace (re^{\pi i x},w), 0 \leq x \leq \frac{2 \pi}{q},|w| \leq \sqrt{2}/2, r \in \R{} \rbrace$  and $ \lbrace (z, re^{\pi i y}), 0 \leq y \leq \frac{2 \pi}{q},|z| \leq \sqrt{2}/2, r \in \R{} \rbrace$.  Each is a regular neighborhood of an arc on the $z$ or $w$-axis and we will refer to these wedges as the $z$-wedge and the $w$-wedge respectively.  The $z$-wedge is pictured in figure \ref{fig:zwedgelink}.  We will show that each of these wedges will glue up to one of the solid tori pictured in figure \ref{knotlens}, and that the arcs of $\D{p/q}$ in these wedges will map to the the arcs of $\LL$ in the solid tori. There are three arcs of $\D{p/q}$ in each of these wedges.  In the $z$-wedge these occur at heights $0$, $\pi/q$ and $2  \pi/q$, where at height $x$, $z= e^{\pi i x}$. The $w$-wedge  also contains three arcs of $\D{p/q}$.  These occur at the levels $0$, $\pi/q$ and $2 \pi/q$, where  $w= e^{\pi i x}$ at level $x$.

\begin{figure}[h]
\centerline{\epsfig{file=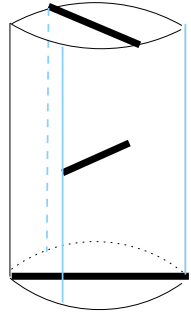, width=1 in}}
\caption{\label{fig:zwedgelink} The link $\D{p/q}$ in the $z$-wedge}
\end{figure}

 Suppose that $q$ is odd.  Consider $(S^3, \D{p/q})$.  As mentioned above, if there is a component of $\D{p/q}$ that intersects the $z$-axis in $(e^{ i x},0)$, then the component that intersects the $z$-axis at $(e^{ i (x + \pi/q)},0)$ is rotated by $p \pi /q$ in the $w$ direction.  Therefore, as in figure \ref{fig:zwedgelink} below, the arc at height $\pi/q$ is rotated by $p \pi/q$ and the arc at height $2 \pi /q$ is rotated by $2p \pi /q$.  This wedge will glue up to a solid torus under the action of $\phi_{p/q}$.  The bottom level will be glued to the top with a twist of $2p \pi /{q}$, and the arc at the bottom will be identified with the arc at the top.  $q$ evenly spaced vertical lines on the boundary of the $z$-wedge will glue up to a $p/q$ curve on the resulting solid torus.   The resulting solid torus is pictured at the top right in figure \ref{knotlens}.    Consider the $w$-wedge of the fundamental domain.  Since $q$ is odd and $(2p,q) =1$, there are integers $x$ and $y$ such that $x (2 \pi p/q) = \pi/q + y \pi$, and $y$ is odd.  The arc at level $\pi / q$ is rotated by $x \pi/ q$, and the arc at level $2 \pi / q$ is rotated by $2x \pi / q$.  Under the action of $\phi_{p/q}$, the bottom and top marked arcs of the $w$-wedge are identified, and the resulting solid torus is pictured in the left in figure \ref{knotlens}.   A meridian curve on this torus will glue to a $p/q$ curve on the solid torus coming from the $z$-wedge.  Therefore we see that the fundamental domain  of the action of $\phi_{p/q}$  on $(S^3, \D{p/q})$ covers $(L_{p/q}, \LL)$, and since $\D{p/q}$ is invariant under this action, $(S^3, \D{p/q})$ covers $(L_{p/q}, \LL)$.

Now let $q$ be even. The $z$-wedge contains three arcs of $D_{p/q}$. The arc at level $\pi/q$ is rotated by $p \pi/q$.  This is $g_{\pi/q, \, p \pi/q}$.  The arc at level $2 \pi/q$ is  part of the  image of the real great circle under $\phi_{p/q}$, and is rotated by $2 \pi p/q$.  Therefore, the $z$-wedge glued up with a $2 \pi p/q$ twist is exactly the solid torus on the right in figure \ref{knotlens}.  The $w$-wedge also glues up to become a solid torus.  Let $n$ be such that $n p \pi  /q =  \pi /q \mod  \pi$.  Then since $q$ is even, $p$ and $n$ must be odd.    The arc at level $2 \pi /q$ is the $n$th in the orbit of the real great circle under $\phi_{p/q}$ and is twisted by $2 \pi n/q$.  The arc at level $\pi/q$ is the $((n-1)/2)$th  in the orbit of $g_{\pi/q, p \pi/q}$.  This is because 
$$\frac{p\pi}{q} + \frac{(n-1)}{2} \frac{ 2p \pi}{q} = \frac{n p \pi}{q} = \frac{\pi}{q} \mod \pi.$$

The arc at level $\pi/q$ in the $w$-wedge is twisted by $\pi/q + (n-1) \pi /q = \pi n /q$.  Since the  top and bottom of the $w$-wedge are identified with a  twist of  $2 \pi n /p$, the $w$-wedge glues up to the solid torus in the left in figure \ref{knotlens}.  A meridian curve of the solid torus coming from the $w$-wedge will be identified with a $p/q$ curve on the boundary of the solid torus coming from the $z$-wedge.  Again, we have that $(S^3, \D{p/q})$ covers $(L_{p/q}, \LL)$. Thus the great circle link complement $S^3 - \D{p/q}$ covers the complement of the knot in the lens space in figure \ref{knotlens}.

Therefore, since $S^3 - \K{p/q}$ is covered by $L_{p/q} -\LL$ and this is in turn covered by the great circle link complement $S^3 - \D{p/q}$, every two-bridge knot and non-trivial link complement is covered by a great circle link complement. By Theorem \ref{fibered}, this finishes the proof of Theorem \ref{twobridgevirtual}.

\end{proof}

By work of Gabai in \cite{detectingfiber}, a rational knot or link $\K{p/q}$ is fibered exactly when $p/q$ has a continued fraction decomposition $ 1/(\pm 2 + 1/( \pm 2 + 1/(\pm 2 +1/(\pm 2 ....)))) $.  Also, the only non-hyperbolic rational knots are the torus knots, (\cite{hatcherthurston}).  Therefore there are infinitely many non-fibered hyperbolic two-bridge knots. In some sense most two-bridge knots are non-fibered.  Thus Theorem \ref{twobridgevirtual} gives a large class of non-fibered, virtually fibered hyperbolic knot complements. 

\section{Some Montesinos knot and link complements \label{montesinossection}} 
\begin{figure}[h]
\centerline{\epsfig{file=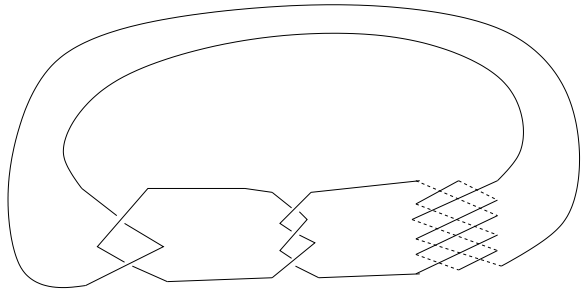, width=3 in}}
\caption{\label{montesinos} A Montesinos knot}
\end{figure}

A {\it Montesinos} knot or link is one that can be written as the union of rational tangles arranged in a circle, as in figure \ref{montesinos}.  The two-fold branched cover of $S^3$ branched along a Montesinos knot or link is a Seifert fibered space.    A proof of this can be found in \cite{Burde}.   A three-manifold is called {\it spherical} if it is a finite quotient of $S^3$ by isometries. We call a Montesinos knot or link  {\it spherical Montesinos} if the associated  Seifert fibered space is spherical. This happens exactly when the base orbifold is a spherical orbifold, and the Euler number is not 0. \cite{scottsurvey, orbifoldbook}. The spherical two-orbifolds with three singular fibers can be computed using the orbifold Euler characteristic \cite{thurstonnotes, orbifoldbook} and are  $\lbrace S(2,2,n), S(2,3,3), S(2,3,4), S(2,3,5) \rbrace$.  For example, the knot in figure \ref{montesinos} is a spherical Montesinos knot, composed of the rational tangles $(1/2)$, $(1/3)$ and $(2/5)$.   

\begin{theorem} 
\label{montesinosfiber}
The complement of any spherical Montesinos knot or link is virtually fibered. \end{theorem}
 \begin{proof}
 
 As in the proof of Theorem \ref{twobridgevirtual}, it suffices to show that any spherical Montesinos  knot or link complement is finitely covered by a great circle link complement.  To this end, let $\K{m}$ be a spherical Montesinos knot or link and $M_k$ denote the double cover of $S^3$ branched along $\K{m}$, with branched covering map $p_1: M_k \rightarrow S^3$.  Then the preimage of the branching locus is a knot or link in $M_k$  which we denote by $\LM$.  Since $M_k$ is a spherical manifold, it is finitely covered by $S^3$ with covering map $p_2: S^3 \rightarrow M_k$.  Now the preimage of $\LM$ is a link in $S^3$ which we denote as $\N$.  We claim that $\N$ is a great circle link.  The composition $p_2 \circ p_1$ is  a branched covering of $S^3$ by $S^3$, where the branching locus is $\K{m}$ and the pre-image of the branching locus in the cover is $\N$.  Call the associated covering group $G$.  Then each component $c$ of $\N$  has the property that the set of elements in $G$ that fix $c$ is cyclic of order two.  The following  special case of W. P. Thurston's orbifold theorem, outlined in  \cite{orbifoldbook}, now applies. 
 
 \begin{theorem}[\cite{orbifoldbook}]
 Assume that $M$ is a orientable irreducible closed three-manifold that contains no incompressible tori.  Suppose that $M$ admits an action by a finite group $G$ of orientation preserving diffeomorphisms such that some non-trivial element has a fixed point set of dimension 1.  Then $M$ has a geometric structure such that this action of the group $G$ is by isometries.  In particular, the fixed point set of each group element is totally geodesic.  
 \end{theorem}
 
 Therefore, the link $\N$ can be realized as a great circle link in the standard metric on $S^3$.  Note that $S^3 -\N$ finitely covers $M_k -\LM$ which two fold covers $S^3 -\K{m}$.  Therefore, $S^3 - \N$ is a finite cover of $S^3 -\K{m}$ which implies by Theorem \ref{fibered} that $S^3 -\K{m}$ is either fibered or virtually fibered.  
 
  \end{proof}

{\bf Remark:} Since the two-fold cover of $S^3$, branched along a two-bridge knot or non-trivial link is also a spherical manifold, we could have also proved Theorem \ref{twobridgevirtual} using the symmetry theorem.  However, this is as not as direct a proof.  Also, it should be possible to prove Theorem \ref{montesinosfiber} using the methods of the proof of Theorem \ref{twobridgevirtual}. The great circle links whose complements cover spherical  Montesinos knot complements are determined in terms of the Grassmannian of two-planes through the origin in \cite{Sakuma}. This paper was only recently discovered by the author. 

 \section{Further Questions  \label{questions}}

We define a knot $k$ to be {\it great} if $S^3 -k$ is commensurable with a great circle link complement. Since virtual fibration is preserved by commensurability, knots that are great are virtually fibered.  In light of Conjecture \ref{virtfiberconj} and Theorem \ref{fibered} the most natural question is:

 \begin{question} \label{everything} Which knots are great? \end{question}

The results of this paper show that all two-bridge knots and all spherical Montesinos knots are great.  The existence of non-great knots is unknown.   However, as a partial answer in this direction we have: 
\begin{theorem} 

There is a knot complement in $S^3$ that is not commensurable with the complement of any strongly invertible great circle link.
\end{theorem}

\begin{proof} 
Let $\Gamma$ be a discrete subgroup of  $\textrm{PSL}_2 \mathbb{C}$ with finite co-volume.  The commensurator of  $\Gamma$ is
$$ Comm (\Gamma) = \lbrace g \in Isom (\mathbb{H}^3): [\Gamma: \Gamma \cap g^{-1}\Gamma g] < \infty \rbrace$$
and $Comm^+(\Gamma)$ is its orientation preserving subgroup. 
By Margulis, $\Gamma$ has a unique maximal element in its commensurability class, $Comm^+(\Gamma) $, if and only if $\Gamma$ is non-arithmetic.  Now let $\Gamma$ be a discrete faithful representation of $\pi_1(S^3 - 9_{32})$.  By Reid, \cite{Reidknot} this group is non-arithmetic.  We will show it is the maximal element in its commensurability class.

In \cite{riley}, Riley shows that the knot $9_{32}$ is asymmetric, meaning that every auto-homeomorphism of the complement  $S^3 - 9_{32}$ is isotopic to the identity. This is done by computing a discrete faithful representation $\Gamma$ of the fundamental group of the knot complement, and showing that there is a fundamental domain that does not admit any hyperbolic symmetries.  Therefore, $\Gamma$ is equal to its normalizer in $\textrm{PSL}_2 \mathbb{C}$.

Therefore, if $\Gamma$ is not the maximal element in its commensurability class, the commensurator must be larger than its normalizer.  Geometrically this means that the knot complement has a hidden symmetry, a symmetry of some finite cover that is not a lift of a symmetry of $S^3 - 9_{32}$.  
Riley shows that the invariant trace field has degree 29, i.e. $[\mathbb{Q}(\textrm{tr} \Gamma): \mathbb{Q}] = 29$.  Since we have a knot complement in $S^3$, the trace field is the same as the invariant trace field, and is an invariant of the commensurability class of $\Gamma$, by Neumann and Reid in \cite{NRarith}.  They also show  \cite[Thm. 9.1]{NRarith} that  a knot complement other than the figure-eight knot complement has hidden symmetries only if the cusp parameter is in $\mathbb{Q}(\sqrt{-1})$ or  $\mathbb{Q}(\sqrt{-3})$.  Since the cusp field is a subfield of the invariant trace field, \cite[Prop. 2.7]{NRarith}, and degree is multiplicative, $S^3 - 9_{32}$ does not have hidden symmetries. 

If $S^3-9_{32}$ was commensurable with a strongly invertible link complement, there would be a finite cover of $S^3 - 9_{32}$  that covered an orbifold.  This orbifold cannot cover $S^3 - 9_{32}$.  Therefore this orbifold either results from a symmetry of $S^3 - 9_{32}$, or a hidden symmetry of $S^3 -9_{32}$.  We have shown that neither of these can happen.

\end{proof}
{\bf Remark:}  It is possible that $S^3 - 9_{32}$ is commensurable with a great circle link that is not strongly invertible.  All great circle links mentioned in this paper are strongly invertible. However, there are great circle links which are not strongly invertible. This is shown in the course of the classification of configurations of lines in $\mathbb{RP}^3$ in \cite{Maz}.  
\vskip .2 in

The fiber exhibited in Theorem \ref{fibered} depended on the choice of a component of a great circle links.  There is one such fiber for each component of the link.  Each of these fibers correspond to an element of $H_2(S^3 - L)$, where $L$ is the great circle link.  A natural question is whether other directions in $H_2(S^3 -L)$ are fibered. In \cite{norm}, W. P. Thurston shows that there is a norm on the second homology of a hyperbolic three-manifold.  Furthermore, he shows that the fibered homology classes are represented by the union of rational lattice points in the cone on some collection of open faces of the unit ball in this norm.  These faces are commonly referred to as the {\it fibered faces}. 

\begin{question} Which faces of the unit ball in the norm on homology of a great circle link complement are fibered faces? 
\end{question} 

{\bf Acknowledgments:} This work began with the author's PhD dissertation, directed by William Thurston.  The author thanks him for his insight, guidance and patience.   More recently, conversations with Alan Reid have  been very helpful. The author also thanks him for a careful reading of the first version of this paper.  This research was supported  by NSF VIGRE grants at the University of California at Davis and the University of Texas at Austin. 

\bibliography{virtuallyfibered1.bib}
 \bibliographystyle{alpha}

 \end{document}